\documentclass[12pt]{amsart}

\newcommand{\E}{{\mathbb E}}
\newcommand{\CI}{{\mathcal I}}

\newcommand{\CZ}{{\mathcal Z}}
\newcommand{\norm}[1]{\lVert #1\rVert}
\newcommand{\nnorm}[1]{\lvert\!|\!| #1|\!|\!\rvert}

\newtheorem{theorem}{Theorem}[section]

\theoremstyle{definition}
\newtheorem{definition}[theorem]{Definition}

\newtheorem{proposition}[theorem]{Proposition}
\newtheorem{corollary}[theorem]{Corollary}
\theoremstyle{remark}

\begin{document}

\title[Averages along cubes for commuting transformations]{Convergence of multiple
ergodic averages along cubes for several commuting transformations}

\author{Qing Chu}

\address{Universit\'{e} Paris-Est, Laboratoire d'Analyse et de math\'{e}matiques
appliqu\'{e}es, UMR
                                               Â´
CNRS 8050, 5 bd Descartes, 77454 Marne la Vall\'{e}e Cedex 2, France}

\email{qing.chu@univ-mlv.fr}


\date{\today}%


\begin{abstract}
We prove the norm convergence of multiple ergodic
averages along cubes for several commuting transformations, and
derive corresponding combinatorial results. The method we use
relies primarily on the ``magic extension" established
recently by B. ~Host.
\end{abstract}

\maketitle

\section{Introduction}
\subsection{Results}
By a system, we mean a probability space endowed with a single or several commuting measure preserving transformations. We prove the following result regarding the convergence of multiple ergodic averages
along cubes for several commuting transformations :

\begin{theorem}\label{main}
Let $d\geq 1$ be an integer and $(X, \mathcal{B}, \mu, T_1, \cdots, T_d)$ be a system. Let $f_{\epsilon},
\epsilon\in\{0,1\}^d\setminus\{00\cdots0\}$ be $2^d-1$ bounded measurable functions on $X$. Then the averages
\begin{align}\label{main1}
\prod_{i=1}^d \frac{1}{N_i-M_i}\sum_{\substack{n_i\in[M_i, N_i)\\
i=1,\dots,d}}\prod_{\substack{\epsilon\in \{0,1\}^d\\
\epsilon\neq 00\cdots0}}T_1^{n_1 \epsilon_1}\cdots
T_d^{n_d\epsilon_d}f_{\epsilon}
\end{align}
converge in $L^{2}(\mu)$ for all sequences of intervals
$[M_1,N_1),\dots,[M_d,N_d)$ whose lengths $N_i-M_i$ ($1\leq i\leq
d$) tend to $\infty$.
\end{theorem}

To illustrate, when $d=2$, the average
(\ref{main1}) is
\begin{equation}\label{dis2}
 \frac{1}{(N_1-M_1)\times(N_2-M_2)}\sum_{\substack{n_1\in[M_1, N_1)\\
n_2\in[M_2, N_2)}}T_1^{n_1}f_{10}\cdot T_2^{n_2}f_{01}\cdot
T_1^{n_1}T_2^{n_2}f_{11}.
\end{equation}

When Theorem ~\ref{main} is restricted to the case that each function $f_{\epsilon}$ is the indicator function of a
measurable set, we have the following lower bound for these
averages:

\begin{theorem}\label{indica}
Let $(X, \mathcal{B}, \mu, T_1, \cdots, T_d)$ be a system and let $A\in
\mathcal{B}$. Then the limit of the averages
\begin{align}\label{indica1}
\prod_{i=1}^d \frac{1}{N_i-M_i}\sum_{\substack{n_i\in[M_i, N_i)\\
i=1,\dots,d}}\mu\big(\bigcap_{\epsilon\in\{0,1\}^d} T_1^{-n_1\epsilon_1}\cdots
T_d^{-n_d\epsilon_d} A\big)
\end{align}
exists and is greater than or equal to $\mu(A)^{2^d}$ for all sequences of intervals
$[M_1,N_1),\dots,[M_d,N_d)$ whose lengths $N_i-M_i$ ($1\leq i\leq
d$) tend to $\infty$.
\end{theorem}

Recall that the upper density $d^*(A)$ of a set $A\subset \mathbb{Z}^d$
is defined to be
$$d^*(A)=\limsup_{\substack{N_i\rightarrow \infty\\1\leq i\leq d}}\prod_{i=1}^d
\frac{1}{N_i}|A\cap[1,N_1]\times\cdots\times [1,N_d]|.$$

A subset $E$ of $\mathbb Z^d$ is said to be syndetic if $\mathbb
Z^d$ can be covered by finitely many translates of $E$.

We have the following corresponding combinatorial result:

\begin{theorem}\label{comb}
Let $A \subset \mathbb{Z}^{d}$ with $d^*(A) > \delta> 0$. Then the set of
$n=(n_1,\dots,n_d)\in \mathbb{Z}^{d}$ such that
$$d^*\left(\bigcap_{\epsilon\in\{0,1\}^d}\{A+(n_1\epsilon_1 ,\dots,
n_d\epsilon_d)\}\right)\geq \delta^{2^d}$$ is syndetic.
\end{theorem}

\subsection{History of the problem}

In the case where $T_1=T_2=\cdots=T_d=T$, the average
(\ref{main1}) is
\begin{equation}\label{diseq}
\prod_{i=1}^d \frac{1}{N_i-M_i}\sum_{\substack{n_i\in[M_i, N_i)\\
i=1,\dots,d}}\prod_{\substack{\epsilon\in \{0,1\}^d\\
\epsilon\neq 00\cdots0}}T^{n_1 \epsilon_1+\cdots
+n_d\epsilon_d}f_{\epsilon}.
\end{equation}

The norm convergence of (\ref{diseq}) was proved by Bergelson for
$d=2$ in ~\cite{B}, and more generally, by Host and Kra for $d>2$ in
~\cite{HK}. The related pointwise convergence problem was studied by Assani and he showed that the averages (\ref{diseq}) converge
a.e. in ~\cite{As1}.

The lower bound for the average (\ref{indica1}) was firstly studied
by Leibman, he provided some lower bounds for (\ref{indica1}) in
~\cite{L}. In the same paper, he gave an example showing that the
average (\ref{dis2}) can diverge if the transformations do not
commute.

However, Assani showed in ~\cite{As1} that the averages
$$\frac{1}{N^2}\sum_{n,m=1}^{N}f(T_1^nx)g(T_2^mx)h(T_3^{n+m}x)$$ do
converge a.e. even if the transformations do not necessarily
commute. He extended this result to the case of six functions in
~\cite{As2}.

The norm convergence of multiple ergodic averages with several
commuting transformations of the form
\begin{equation}\label{retao}
\frac{1}{N}\sum_{n=1}^N T_1^n f_1\cdot\ldots\cdot T_d^n f_d,
\end{equation}
was proved by Conze and Lesigne ~\cite{CL} when $d=2$. The general
case was originally proved by Tao ~\cite{H}, and subsequent proofs
were given by Austin ~\cite{A}, Host ~\cite{H} and Towsner~
\cite{Tow}.

\subsection{Methods}
The main tools we use in this paper are the seminorms and the existence
of ``magic extensions" for commuting transformations established by
Host ~\cite{H}. The ``magic extensions" can be viewed as a concrete
form of the pleasant extensions built by Austin in ~\cite{A}.

\subsection*{Acknowledgement}
This paper was written while the author was visiting the Mathematical Sciences Research Institute and the author is grateful for their kind hospitality. The author also thanks the referee of this paper for remarks.

\section{Seminorm and upper bound}
\subsection{Notation and definitions}\label{nota}
For an integer $d\geq 1$, we write $[d]=\{1,2,\dots,d\}$ and identify $\{0,1\}^d$ with the
family of subsets of $[d]$. Therefore, the assertion ``$i\in
\epsilon$'' is equivalent to $\epsilon_i=1$. In particular,
$\emptyset$ is the same as $00\cdots0\in\{0,1\}^d$. We write
$|\epsilon|=\sum_i \epsilon_i$ for the number of elements in
$\epsilon$.

Let $(X, \mu, T_1, \cdots, T_d)$ be a system. For each $n=(n_1,\dots,n_d)$, $\epsilon=\{i_1,\dots,i_k\}\subset
[d]$, and for each integer $1\leq k\leq d$, we write
$$T_{\epsilon} ^n= T_{i_1}^{n_{i_1}}\cdots T_{i_k}^{n_{i_k}}.$$

For any transformation $S$ of some probability space, we denote by $\CI {(S)}$ the
$\sigma$-algebra of $S$-invariant sets.

We define a measure $\mu_1$ on $X^2$ by
$$\mu_1=\mu\times_{\CI (T_1)}\mu_1.$$

This means that for $f_0,f_1\in L^{\infty}(\mu)$, we have
$$\int (f_0\otimes f_1) (x_0,x_1) d\mu_1(x_0,x_1)=\int \E(f_0|\CI(T_1))\cdot\E(f_1|\CI(T_1))d\mu.$$

For $2\leq k\leq d$, we define a measure $\mu_{k}$ (see \cite{H}) on
$X^{2^{k}}$ by
$$\mu_k=\mu_{k-1}\times_{\CI (T^{\triangle}_k)}\mu_{k-1},$$
where $T^{\triangle}_k :=\underbrace{T_k\times\cdots\times
T_k}_{2^{k-1}}.$

We write $X^*=X^{2^d}$, and points of $X^*$ are written as
$x=(x_\epsilon:\epsilon\subset [d])$. We write $\mu^*:=\mu_d$.

For $f\in L^{\infty}(\mu)$, define
$$\nnorm{f}_{T_1,\dots,T_d}:=\left(\int \prod_{\epsilon\in\{0,1\}^d}
f(x_{\epsilon})d\mu^*(x)\right)^{1/2^{d}}.$$

It was shown in Proposition 2 in ~\cite{H} that
$\nnorm{\cdot}_{T_1,\dots,T_d}$ is a seminorm on $L^{\infty}(\mu)$.
We call this the box seminorm associated to $T_1,\dots,T_d.$

For $\epsilon\subset[d]$, $\epsilon\neq\emptyset$, we write
$\nnorm\cdot_\epsilon$ for the seminorm on $L^\infty(\mu)$
associated to the transformations $T_i$, $i\in\epsilon$. For
example, $\nnorm\cdot_{110\cdots 00}$ is  the seminorm associated to
 $T_1,T_2$ because $\epsilon=110\cdots0\in\{0,1\}^d$ is identified with
$\{1,2\}\subset[d]$.

By Proposition 3 in ~\cite{H}, if we rearrange the order of the
digits in $\epsilon$, the seminorm $\nnorm \cdot_{\epsilon}$ remains
unchanged.

\subsection{Upper bound}
In the following, we assume that all functions $f_\epsilon$,
$\epsilon\subset[d]$, are real valued and satisfy $|f_\epsilon|\leq
1$.

\begin{proposition}\label{uppba} Maintaining the above notation and hypotheses,
\begin{align}\label{uppba1}
 \limsup_{\substack{N_i-M_i\rightarrow\infty\\i=1,\dots,d}}\Bigl\Vert\prod_{i=1}^d
\frac{1}{N_i-M_i}
 \sum_{n\in [M_1,N_1)\times\cdots\times [M_d,N_d)}\prod_{\substack{\epsilon\subset
[d]\\
\epsilon\neq \emptyset}}T_{\epsilon}^{n}f_{\epsilon}\Bigr\Vert_{L^2(\mu)}\leq
\min_{\substack{\epsilon\subset [d]\\
\epsilon\neq \emptyset}} \nnorm {f_{\epsilon}}_{T_1\dots T_d}.
\end{align}
\end{proposition}

\begin{proof}
We proceed by induction on $d$. For $d=1$, we have
$$\Bigl\Vert\frac{1}{N_1-M_1}\sum_{n_1\in [M_1,N_1)}T_1^{n_1}
f_1\Bigr\Vert^2_{L^2(\mu)}\rightarrow \int
\E(f_1|\mathcal{I}(T_1))^2 d\mu= \nnorm {f_1}^2_{T_1}.$$

Let $d\geq 2$ and assume that (\ref{uppba1}) is established for
$d-1$ transformations.

We show that for every $\alpha\subset [d], \alpha\neq
\emptyset$, the limsup of the the left hand side of (\ref{uppba1})
is bounded by $\nnorm{f_{\alpha}}_{T_1,\dots,T_d}$. By a permutation of digits if needed, 
we can assume that $\alpha\neq \underbrace{0\dots
0}_{d-1}1$. The square of the norm in the left hand side of
(\ref{uppba1}) is equal to

\begin{align*}
\Bigl\Vert\frac{1}{N_d-M_d}\sum_{n_d\in
[M_d,N_d)}T_{d}^{n_d}f_{0...01}\cdot\prod_{i=1}^{d-1}
\frac{1}{N_i-M_i}\sum_{m\in [M_1,N_1)\times\cdots\times
[M_{d-1},N_{d-1})}\\\prod_{\substack{\eta\subset [d-1]\\
\eta\neq\emptyset}}T^{m}_{\eta}(f_{\eta 0}\cdot T_d^{n_d}f_{\eta
1})\Bigr\Vert_{L^{2}(\mu)}^{2}.
\end{align*}

By the Cauchy-Schwartz Inequality, this is less than or equal to
\begin{align}\label{uppba3}
\begin{split}
\frac{1}{N_d-M_d}\sum_{n_d\in
[M_d,N_d)}\Bigl\Vert\prod_{i=1}^{d-1}
\frac{1}{N_i-M_i}\sum_{m\in [M_1,N_1)\times\cdots\times
[M_{d-1},N_{d-1})}\\\prod_{\substack{\eta\subset [d-1]\\
\eta\neq\emptyset}}T^{m}_{\eta}(f_{\eta 0}\cdot T_d^{n_d}f_{\eta
1})\Bigr\Vert_{L^{2}(\mu)}^{2}.
\end{split}
\end{align}

By the induction hypothesis, when $N_i-M_i\rightarrow\infty$,
$i=1,\dots,d-1$, the limsup of the square of the norm in (\ref{uppba3}) is less than or equal to
$$\min_{\substack{\eta\subset [d-1]\\
\eta\neq \emptyset}}\|f_{\eta 0}\cdot T_d^{n_d}f_{\eta
1}\|^2_{T_1,\dots,T_{d-1}},
$$
where $\|\cdot\|_{T_1,\dots,T_{d-1}}$ is the seminorm associated to
the $d-1$ transformations $T_1,\dots, T_{d-1}$.

Note that $\alpha$ is equal to $\eta 0$ or $\eta 1$ for some
$\eta\subset [d-1]$, and by the Cauchy-Schwartz Inequality, we have
\begin{align*}
 &\ \
\lim_{N_d-M_d\rightarrow\infty}\frac{1}{N_d-M_d}\sum_{n_d\in[M_d,N_d)}\|f_{\eta
0}\cdot T_d^{n_d}
 f_{\eta
1}\|^{2^{d-1}}_{T_1,\dots,T_{d-1}}\\&=\lim_{N_d-M_d\rightarrow\infty}\frac{1}{N_d-M_d}\sum_{n_d\in[M_d,N_d)}
 \int \bigotimes_{\eta\subset [d-1]}(f_{\eta 0}\cdot T_d^{n_d}f_{\eta 1}) d\mu_{d-1}\\
 &=\int \mathbb{E}(\bigotimes_{\eta\subset [d-1]} f_{\eta
0}|\mathcal{I}(T^{\triangle}_d)\cdot\mathbb{E}
 (\bigotimes_{\eta\subset [d-1]} f_{\eta 1}|\mathcal{I}(T^{\triangle}_d)d\mu_{d-1} \\
 &\leq\left(\int
 |\mathbb{E}(\bigotimes_{\alpha\subset
[d]}f_{\alpha}|\mathcal{I}(T^{\triangle}_{d}))|^{2}d\mu_{d-1}\right)^{1/2}
 \\&=\left(\int \bigotimes_{\alpha\subset [d]} f_{\alpha}
 d\mu^{*}\right)^{1/2}=\nnorm{f_{\alpha}}^{2^{d-1}}_{T_1,\dots,T_d}.
\end{align*}
This completes the proof.
\end{proof}

The following proposition is a generalization of Proposition
\ref{uppba}, although its proof depends upon Proposition
\ref{uppba}.

\begin{proposition}\label{uppbb}
Let $r$ be an integer with $1\leq r\leq d$. Then
\begin{equation}\label{uppbb1}
 \limsup_{\substack{N_i-M_i\rightarrow\infty\\i=1,\dots,d}}
\Bigl\Vert \prod_{i=1}^d \frac{1}{N_i-M_i} \sum_{n\in
[M_1,N_1)\times\cdots\times [M_d,N_d)}
   \prod_{\substack{\epsilon\subset[d]\\0<|\epsilon|\leq r}}
T_{\epsilon}^{n}f_{\epsilon} \Bigr\Vert_{L^{2}(\mu)} \leq
\min_{\substack{\epsilon\subset[d]\\ |\epsilon|=r}}
\nnorm{f_\epsilon}_{\epsilon}
\end{equation}
\end{proposition}

\begin{proof}
We show that for every $\alpha\subset[d]$ with $|\alpha|=r$, the $\limsup$ in (\ref{uppbb1}) is bounded by
$\nnorm{f_\alpha}_\alpha$ .

By a permutation of digits we can restrict to the case that
$$\alpha=\underbrace{11\dots 1}_{r}00\dots 0\ .$$

We show that

\begin{equation}\label{uppbb2}
 \limsup_{\substack{N_i-M_i\rightarrow\infty\\i=1,\dots,d}}\left\|
 \prod_{i=1}^d \frac{1}{N_i-M_i}\sum_{n\in [M_1,N_1)\times\cdots\times [M_d,N_d)}
 \prod_{\substack{\epsilon\subset [d]\\
0<|\epsilon|\leq r}}T_{\epsilon}^{n}f_{\epsilon}\right\|_{L^{2}(\mu)}\leq
\nnorm{f_{\alpha}}_{\alpha}.
\end{equation}

The norm in (\ref{uppbb2}) is equal to

\begin{align}\label{uppbb3}
\begin{split}
\big\| \prod_{i=r+1}^{d}\frac{1}{N_i-M_i}\sum_{m\in
[M_{r+1},N_{r+1})\times\cdots\times
[M_d,N_d)}\prod_{\substack{\epsilon\subset \{r+1,\dots,d\}\\
\epsilon\neq
\emptyset}}T_{\epsilon}^{m}f_{\epsilon}\cdot\prod_{j=1}^{r}\frac{1}{N_j-M_j}
\\\sum_{n\in [M_1,N_1)\times\cdots\times [M_r,N_r)}\big(\prod_{\substack{\eta\subset
[r]\\
\eta\neq\emptyset}}T^{n}_{\eta} \prod_{\substack{\theta\subset [d-r]\\
|\eta\theta|\leq r}}T^{n}_{r+\theta}f_{\eta \theta}\big)\cdot
(T_{1}^{n_1}\cdots T_{r}^{n_r}f_{\alpha})\big\|_{L^{2}(\mu)}.
\end{split}
\end{align}
where $r+\theta=\{r+k\colon \ k\in\theta\}$.

Let \[g_{\eta}= \begin{cases}\displaystyle\prod_{\substack{\theta\subset[d-r] \\
|\eta\theta|\leq r}}T^{n}_{r+\theta}f_{\eta \theta} & 0<|\eta|<r\ ; \\
f_{\alpha} & |\eta|=r\ .\end{cases}
\]

Then (\ref{uppbb3}) is equal to

\begin{align}\label{uppbb4}
\begin{split}
\Bigl\Vert
 \prod_{i=r+1}^{d}\frac{1}{N_i-M_i}
\cdot \sum_{m\in [M_{r+1},N_{r+1})\times\cdots\times [M_d,N_d)}
\Bigl(\prod_{\substack{\epsilon\subset \{r+1,\dots,d\}\\
\epsilon\neq\emptyset}}
T^{m}_{\epsilon}f_{\epsilon}\Bigr)\\
\cdot\prod_{j=1}^{r}\frac{1}{N_j-M_j} \cdot\sum_{n\in
[M_1,N_1)\times\cdots\times [M_r,N_r)}
\Bigl(\prod_{\substack{\eta\subset [r]\\ \eta\neq\emptyset}}
T^{n}_{\eta} g_{\eta}\Bigr)\Bigr\Vert_{L^{2}(\mu)}.
\end{split}
\end{align}

By the Cauchy-Schwartz Inequality, the square of (\ref{uppbb4}) is less than or equal to
\begin{align}\label{uppbb5}
\prod_{i=r+1}^{d}\frac{1}{N_i-M_i}\sum_{m\in
[M_{r+1},N_{r+1})\times\cdots\times
[M_d,N_d)}\Bigl\Vert\prod_{j=1}^{r}\frac{1}{N_j-M_j}\sum_{n\in
[M_1,N_1)\times\cdots\times [M_r,N_r)}
\prod_{\substack{\eta\subset [r]\\
\eta\neq\emptyset}}T^{n}_{\eta} g_{\eta}\Bigr\Vert_{L^{2}(\mu)}^{2}.
\end{align}

By Proposition \ref{uppba}, the limsup of (\ref{uppbb5}) as
$N_i-M_i\rightarrow\infty$, $i=1,\dots,r$ is bounded by

\begin{align*}
\prod_{i=r+1}^{d}\frac{1}{N_i-M_i}\sum_{\substack{n_i\in[M_i, N_i)\\
i=r+1,\dots,d}}\|f_{\alpha}\|_{T_1,\dots, T_r}^{2}=\|f_{\alpha}\|_{T_1,\dots, T_r}^{2}.
\end{align*}

This completes the proof.
\end{proof}

\section{The case of the magic extension}

We recall the definition of a ``magic'' system.

\begin{definition}[Host, ~\cite{H}]
A system $(X,\mu,T_1,\dots,T_d)$ is called a ``magic'' system if $f\in L^{\infty}(\mu)$ is such that $\E(f|\bigvee_{i=1}^d\CI(T_i))=0$, then $\nnorm f_{T_1,\dots,T_d}=0$.
\end{definition}

Given a system $(X,\mu,T_1,\dots,T_d)$, let $X^*$ and $\mu^*$ be defined as in Section ~\ref{nota}. 
We denote by $T_i^*$ the \emph{side transformations} of $X^*$, given by
\[\text{for every }\epsilon\in\{0,1\}^d,\ \
(T_i^*x)_{\epsilon}=\begin{cases}T_ix_\epsilon & \text{if}\ \epsilon_i=0\ ;\\
x_\epsilon & \text{if}\ \epsilon_i=1\ .\end{cases}\]

By Theorem 2 in ~\cite{H}, $(X^{*},\mu^*, T_1^*,\dots,T_d^*)$ is a ``magic''
system, and admits $(X,\mu,T_1,\dots,T_d)$ as a factor.


For  $\epsilon\subset[d]$, $\epsilon\neq\emptyset$, we write $\nnorm
\cdot^*_\epsilon$ for the seminorm on $L^\infty(\mu^*)$ associated
to the transformations $T^*_i$, $i\in\epsilon$. Moreover, we define the
$\sigma$-algebra
$$
 \CZ^*_\epsilon:=\bigvee_{i\in\epsilon}\CI(T_i^*)
$$
of $(X^{*},\mu^{*})$. For example, $\CZ^*_{\{1,2,d\}}=\CI(T_1^*)\vee\CI(T_2^*)\vee
\CI(T_d^*)$.

We prove Theorem ~\ref{main} for the magic system $(X^*,\mu^*,
T_1^*,\dots,T_d^*)$. 

\begin{theorem}\label{mainm}
Let $f_\epsilon$, $\epsilon\subset[d]$, be functions on $X^*$
with $\norm{f_\epsilon}_{L^\infty(\mu^*)}\leq 1$ for every
$\epsilon$. Then the averages
\begin{equation}\label{mainm1}
\prod_{i=1}^d\frac 1{N_i-M_i} \sum_{n\in [M_1,N_1)\times\cdots\times
[M_d,N_d)}
\prod_{\substack{\epsilon\subset [d]\\
\epsilon\neq \emptyset}}T^{*n}_{\epsilon}f_{\epsilon}
\end{equation}
converge in $L^2(\mu^*)$ for all sequences of intervals
$[M_1,N_1),\dots,[M_d,N_d)$ whose lengths $N_i-M_i$ ($1\leq i\leq
d$) tend to $\infty$.
\end{theorem}

Since the system $(X^{*},\mu^*, T_1^*,\dots,T_d^*)$ admits $(X,\mu,T_1,\dots,T_d)$ as 
a factor, Theorem ~\ref{mainm} implies our main result Theorem ~\ref{main}.

\begin{theorem}\label{charac}
For every $\epsilon\subset[d]$, $\epsilon\neq\emptyset$, and every
function $f\in L^\infty(\mu^*)$, we have:
\begin{equation}\label{charac1}
\text{If  }\E_{\mu^*}(f\mid\CZ^*_\epsilon)=0,\text{ then } \nnorm
f_\epsilon^*=0.
\end{equation}
\end{theorem}

\begin{proof}
Assume $|\epsilon|=r>0$. By a permutation of digits we can
assume that
$$
 \epsilon=\{ d-r+1,d-r+2,\dots,d\}.
$$

We define a new system $(Y, \nu, S_1,\dots,S_r)$, where
$Y={X}^{2^{d-r}}$, $\nu=\mu_{d-r}$, the $d-r$ step measure
associated to $T_1^*,\dots,T_{d-r}^*$. Define
$$S_i=\underbrace{T_{d-r+i}\times\cdots\times T_{d-r+i}}_{2^{d-r}}$$
on $Y$ for $i=1,\dots,r$.

Note that by definition, $Y^*=X^{2^d}=X^{*}$, and
$$S^*_i=T^*_{d-r+i},\ \ S^{\triangle}_i=T^{\triangle}_{d-r+i}$$
for $i=1,\dots,r.$ Moreover,
$$\nu_1=\nu\times_{\CI (S_1)}\nu=\mu_{d-r}\times_{\CI
(T^{\triangle}_{d-r+1})}\mu_{d-r}=\mu_{d-r+1}.$$ By induction,
$$\nu_{i+1}=\nu_{i}\times_{\CI (S_{i+1}^{\triangle})}\nu_{i}=\mu_{d-r+i}\times_{\CI
(T_{d-r+i+1}^{\triangle})}\mu_{d-r+i}=\mu_{d-r+i+1},$$ for
$i=1,\dots,r-1$.

Therefore $(X^{*}, \mu^*, T^*_{d-r+1},\dots, T^*_{d})$ is just the
magic extension $(Y^*, \nu^*, S^*_1,\dots,S^*_r)$ of
$(Y,\nu,S_1,\dots,S_r)$. So
$$\mathcal Z_{\epsilon}^*=\bigvee_{i\in\epsilon}\CI(T_i^*)=\bigvee_{i=1}^r\CI(S_i^*):=\mathcal{W}^*_Y.$$

If $f\in L^\infty(\mu^*)$ with $\E_{\mu^*}(f\mid\CZ^*_\epsilon)=0$, this
is equivalent to $\E_{\mu^*}(f\mid\mathcal{W}^*_Y)=0$, and by Theorem 2
in \cite{H}, we have $\nnorm {f}^*_{S_1^*,\dots,S_r^*}=0$. Thus
$\nnorm f_\epsilon^*=\nnorm {f}^*_{S_1^*,\dots,S_r^*}=0$
\end{proof}

\begin{proposition}\label{object}
Let $f_\epsilon$, $\epsilon\subset[d]$, be functions on $X^*$
with $\norm{f_\epsilon}_{L^\infty(\mu^*)}\leq 1$ for every
$\epsilon$. Let $r$ be an integer with $1\leq r\leq d$. Then the
averages
\begin{equation}\label{object1}
\prod_{i=1}^d\frac 1{N_i-M_i} \sum_{n\in [M_1,N_1)\times\cdots\times
[M_d,N_d)}
\prod_{\substack{\epsilon\subset [d]\\
0<|\epsilon|\leq r}}T^{*n}_{\epsilon}f_{\epsilon}
\end{equation}
converge in $L^2(\mu^*)$ for all sequences of intervals
$[M_1,N_1),\dots,[M_d,N_d)$ whose lengths $N_i-M_i$ ($1\leq i\leq
d$) tend to $\infty$.
\end{proposition}

We remark that Theorem \ref{mainm} follows immediately from this
proposition when $r=d$.

\begin{proof}
We proceed by induction on $r$. When $r=1$, the average (\ref{object1}) is
\begin{equation}
 \prod_{i=1}^d \frac{1}{N_i-M_i}\sum_{\substack{n_i\in[M_i, N_i)\\
i=1,\dots,d}}T_1^{*n_1}f_{10\dots 0}\cdots T_d^{*n_d}f_{0\dots 01}.
\end{equation}
By the Ergodic Theorem, this converges to $\E(f_{10\dots 0}|\CI(T_1^*))\cdots\E(f_{0\cdots 01}|\CI(T_d^*))$.

Assume $r>1$, and that the proposition is true for $r-1$
transformations.

For $\alpha\subset [d]$, $|\alpha|=r$, if
$\E_{\mu^*}(f_\alpha\mid\CZ^*_\alpha)=0$, then by Theorem
~\ref{charac}, we have $ \nnorm {f_{\alpha}}_{\alpha}^*=0$. By
Proposition ~\ref{uppba}, the average (\ref{object1})
converges to $0$. Otherwise, by a density argument, we can assume that
$$f_\alpha=\prod_{i\in\alpha}f_{\alpha,i}$$ where
$f_{\alpha,i}$ is $T_i^*$-invariant. Then
$$T_{\alpha}^{*n}f_{\alpha}=\prod_{i\in\alpha}
T^{*n}_{\alpha\backslash\{i\}}f_{\alpha,i}\ .$$

Thus
$$\prod_{\substack{\epsilon\subset [d]\\
0<|\epsilon|\leq r}}T^{*n}_{\epsilon}f_{\epsilon}=\prod_{\substack{\eta\subset [d]\\
0<|\eta|\leq {r-1}}}T^{*n}_{\eta}g_{\eta}\ ,$$ where
\[g_{\eta}=\begin{cases}f_{\eta} & |\eta|<r-1\ ; \\f_{\eta}\displaystyle\prod_{i\notin
\eta}f_{\eta\cup i,i} & |\eta|=r-1\ .\end{cases}\]

Therefore (\ref{object1}) converges by the induction hypothesis.
\end{proof}

\section{combinatorial interpretation}
\begin{proof}[Proof of Theorem ~\ref{indica}]
Apply Theorem ~\ref{main} to the indicator function $\textbf{1}_{A}$,
we know that the limit of the averages
\begin{equation}\label{indica2}
\prod_{i=1}^d \frac{1}{N_i-M_i}\sum_{\substack{n_i\in[M_i, N_i)\\
i=1,\dots,d}}\int\prod_{\epsilon\in \{0,1\}^d}T_1^{n_1
\epsilon_1}\cdots T_d^{n_d\epsilon_d}\textbf{1}_{A}d\mu
\end{equation}
exist. By Lemma 1 in ~\cite{H}, if we take the limit as
$N_1-M_1\rightarrow\infty$, then as $N_2-M_2\rightarrow\infty$, ...
and then as $N_d-M_d\rightarrow\infty$, the average (\ref{indica2})
converges to $\nnorm {\textbf{1}_{A}}^{2^d}_{T_1,\dots,T_d}$. Thus the
limit of the average (\ref{indica2}) is $\nnorm
{\textbf{1}_{A}}^{2^d}_{T_1,\dots,T_d}$. Since

$$\nnorm {f}^{2^d}_{T_1,\dots,T_d}=\Vert\E (\bigotimes_{\epsilon\subset [d-1]}f
|\CI (T_d^{\triangle}))\Vert_{L^2(\mu_{d-1})}^2\geq(\int
\bigotimes_{\epsilon\subset [d-1]}f d\mu_{d-1})^2=\nnorm
{f}^{2^{d}}_{T_1,\dots,T_{d-1}},$$ we have $\nnorm
{\textbf{1}_{A}}_{T_1,\dots,T_d}\geq \nnorm
{\textbf{1}_{A}}_{T_1}\geq\int \textbf{1}_{A}d\mu=\mu(A)$, and the
result follows.
\end{proof}

Theorem \ref{indica} has the following corollary:
\begin{corollary}\label{indicacor}
Let $(X, \mathcal{B}, \mu, T_1,\dots,T_d)$ be a system, where $T_1,\dots, T_d$
are commuting measure preserving transformations, and let $A
\in\mathcal{B}$. Then for any $c > $0, the set of $n\in\mathbb{Z}^k$
such that
$$\mu\big(\bigcap_{\epsilon\in\{0,1\}^d} T_1^{-n_1\epsilon_1}\cdots
T_d^{-n_d\epsilon_d} A\big)\geq \mu(A)^{2^d}-c$$
is syndetic.
\end{corollary}
The proof is exactly the same as Corollary 13.8 in ~\cite{HK}.

Theorem ~\ref{comb} follows by combining Furstenberg's correspondence
principle and Corollary ~\ref{indicacor}.

\bibliographystyle{amsplain}

\end{document}